\documentclass[a4paper,12pt]{amsart}
\usepackage{latexsym}
\usepackage[french,english]{babel}
\usepackage{epic,eepic,eucal}
\usepackage[T1]{fontenc}          
\usepackage{enumerate}

\def\s{\mathcal S}

\def\a{\alpha}

\def\R{\mathop{\bf R}}
\def\Z{\mathop{\bf Z}}
\def\T{\mathop{\bf T}}
\def\N{\mathop{\bf N}}

\def\X{ \mathcal L}

\def\A{          \mathcal A}
\def\d{ \delta}

\def\E{          \mathcal E}

\def\A{          \mathcal A}

\def\a{         \alpha}
\def\s{\mathcal S}

\def\R{\mathop{\bf R}}
\def\Q{\mathop{\bf Q}}

\def\Z{\mathop{\bf Z}}

\def\carre{ \hfill $\Box$}

\usepackage{enumerate}
\makeindex

\renewcommand{\subsubsection}{\vspace{0.3cm}
  \refstepcounter{subsection} \noindent   
  \thesubsubsection \ }

\newtheorem{defi}{\sc Definition.}

\begin{document}
\title{Continuous spectrum on laminations over 
Aubry-Mather sets.} \label{trois1}
\author{Bassam R. Fayad}
\address{Bassam  Fayad, LAGA, Universit\'e Paris 13, Villetaneuse, }
\email{fayadb@math.univ-paris13.fr}

\maketitle
\begin{abstract} {If we perturb a completely integrable
  Hamiltonian system with two degrees of freedom, the perturbed flow    
  might display, on every energy level, invariant sets that        
  are laminations  over Aubry-Mather sets of a Poincar{\'e} section of the flow. Each one of these laminations  carry a unique invariant probability measure for the flow and it is interesting therefore to understand the  statistical properties of these measures.  From a result of Kocergin in \cite{Ko1}, we know that mixing is {\sl a priori} impossible. In this paper, we investigate on the possible occurrence of weak mixing.

The answer will essentially depend on the number of orbits of gaps in the Aubry-Mather set. More precisely, if the Aubry-Mather set has exactly one orbit of gaps  and is hyperbolic then the special flow over it with any smooth ceiling function will be conjugate to a suspension with a constant ceiling function, failing hence to be weak mixing or even topologically weak mixing.  
  To the contrary,
  if the Aubry-Mather set has more than one orbit of gaps with at least two in a {\sl general
  position} then the special flow over it will in general be weak mixing.}
\end{abstract}

\section{Introduction.}

For an integrable Hamiltonian system with $N$ degrees of freedom, the
$2N$-dimensional phase space is completely foliated with invariant
$N$-dimensional  tori on which  the motion is that of a translation. The
Kolmogorov-Arnold-Moser (KAM) theorem states that for a generic and sufficiently
small $C^{\infty}$
 perturbation of such a completely integrable system,  an
arbitrarily large proportion (in measure) of the invariant tori will
be preserved, i.e. slightly deformed into invariant tori  where  the
dynamics is $C^{\infty}$ conjugate to  translation flows.
These preserved tori correspond in general to frequencies that are badly approximated by rational vectors, i.e.  Diophantine translation vectors.
It is then a natural question to ask what type of behavior can we expect on the other invariant sets of the perturbed system, sets corresponding to the non-Diophantine frequencies. 

Among such sets we find, in the closure of the KAM tori, invariant tori that carry non linearizable flows.  Indeed, for the   generic
$C^{\infty}$ perturbation of the integrable Hamiltonian system,   M. Herman proved in \cite{He3} that 
 on a residual subset among 
the invariant tori of the perturbed flow the dynamics is
 uniquely ergodic and weak mixing. Since the invariant tori he considers lie in the closure of KAM tori, to prove the result above Herman investigates on the generic behavior of a diffeomorphism that is in the closure of smoothly linearizable ones. To this end, one way is to use the successive conjugations techniques as in \cite{AK} or \cite{HeF} which imply genericity of unique ergodicity and weak mixing. Another way, which can still be seen as a particular case of successive conjugations, is to consider reparametrizations of linear flows since such reparametrized flows  are always  in the closure of flows that are smoothly conjugate to linear ones (as it follows from the famous result of Kolmogorov on the rigidity under smooth reparametrizations of Diophantine linear flows \cite{Kkk}). Following this program we have studied in 
 \cite{weakmixing} the occurrence of weak mixing for reparametrizations
of minimal translation flows on the torus ${\T}^n$, $n \geq 2$ and
 proved that for a non-Diophantine linear flow, the flow obtained after reparametrization  with a
strictly positive smooth function $\phi$ is in general  (i.e. for
$\phi$ in a $G_{\delta}$ dense subset of $C^{\infty}( {\T}^n, {\R}_+^*)$) weak mixing for
its unique invariant measure. That is the content of the theorem
 we will  state  in \S \ref{tores}.

\vspace{0.15cm}

What about the dynamics on invariant sets that are not tori?
Here, we will only consider Hamiltonian systems with two degrees of
freedom and a class of perturbations of completely integrable systems for which the perturbed flow displays invariant sets that are laminations
 over invariant Aubry-Mather sets. Indeed,
 consider on an energy level of the perturbed flow a two dimensional Poincar{\'e}
section: by KAM theory, this section contains invariant circles and hence an   invariant annulus for the Poincar\'e return map. In addition, the Poincar\'e map is 
under generic conditions on the perturbation a twist map for 
which the Aubry-Mather theory is thus valid   (Cf. Section \ref{amsetting}). Each Aubry-Mather set carries a unique invariant probability measure for the Poincar\'e map, and so does the laminations  for the suspended flow.  For the perturbations we consider the dynamics on the laminations is not weak mixing because it is a constant time suspension. 
Our  aim in this paper is to study the display  of weak mixing after time change (Cf. \S \ref{rep}  and \ref{weak} for the definitions).

Depending on whether the Aubry-Mather set (we will call
Aubry-Mather sets the Cantor minimal invariant sets exclusively) has one or many {\sl holes}  (Cf. Definition \ref{PPPP} in \S \ref{trous})   the
results will be very contrasting. More
precisely, given an Aubry-Mather set $\A_\a$ of the Poincar{\'e} map, we
  denote by $\s_\a$ the corresponding invariant laminations by the
  flow and  prove the following results: if the Aubry-Mather set $\A_\a$ has one
  hole and is hyperbolic, then any smooth reparametrization of the flow has its restriction to $\s_\a$ conjugate to a suspension flow over $\A_\a$ with  a constant
  ceiling function, failing thus to be weak mixing (Cf. Theorem \ref{amoneone} in \S \ref{amresults}). To the contrary,
  if $\A_\a$ has more than one hole with at least two in a {\sl general  position} (Cf. Definition \ref{amhhhh} in \S \ref{trous}), then arbitrarily close to 1 (in the $C^{\infty}$ topology),
  there exists a function $\phi$ such that the reparametrization of the perturbed  flow by $\phi$ is weak mixing on $\s_\a$ (Theorem \ref{ammanymany}).

\vspace{0.2cm} 

\noindent {\sc Remark 1.}  We do not know how often two holes of an Aubry-Mather set are in a general position. The  arithmetic condition between the point projections of two gaps that defines  them as being in general position (Cf. Definition \ref{amhhhh}) is generic but we do not know how to lift the condition into a generic property on the Aubry-Mather sets.  However, we will see that for any irrational $\a$, 
the numbers $0$ and ${1 / 2}$ are in a general position
with respect to $\a$ (Cf. the proof of Proposition \ref{am3}). 
 Hence, considering two fold coverings, an Aubry-Mather set with one hole becomes an Aubry-Mather with two holes (one in $0$, and one in $1 / 2$) in a general position for which Theorem \ref{ammanymany} applies (Cf. Corollary \ref{co}).

\vspace{0.2cm}

\noindent {\sc Remark 2.} It is worth noticing here that an Aubry-Mather set is typically
hyperbolic (as proved by Patrice Le Calvez \cite{Patrice}, Cf. \S 2.7.1  below) and has one hole (the single gap Theorem by Veerman
\cite{veerman}). Hence the absence of weak mixing prevails on the laminations.
Furthermore, Theorem \ref{amoneone} can be viewed as a statement of  rigidity under time change in the case of non-Diophantine frequencies when the invariant tori degenerate into Aubry-Mather sets, similar to the rigidity result of Kolmogorov in the case of invariant tori carrying Diophantine linear flows \cite{Kkk}, and surprisingly contrasting with the instability results for  invariant tori with Liouvillian frequencies \cite{He3}, \cite{weakmixing} (See \S  \ref{tores}). 



\section{The setting.} \label{amsetting} 

\subsection{} \label{set} Let $(\theta,r,s,u)$ be a system of coordinates on ${\T} \times \R \times
{\T} \times \R$. Given a real function $H(\theta,r,s)$ of class
$C^{\infty}$, and a real number $\varepsilon$, we consider the smooth Hamiltonian system
given by
\begin{eqnarray} \hat{H} (\theta,r,s,u) := u + { r^2 \over 2} + \varepsilon H(\theta,r,s). \label{hamo} \end{eqnarray}

We consider the standard symplectic structure $\omega_0 = d\theta \wedge dr + ds \wedge du$. We denote by $X_{ \hat{H}}^t$ the Hamiltonian
flow corresponding to $\hat{H}$, that is the flow given by the vector field $X_{\hat{H}}$ satisfying the formula $\omega_0 (X_{\hat{H}},.) = d \hat{H} (.)$.  Since $ { ds / dt} = {\partial \hat{H} / \partial u} = 1,$
the coordinate $s$ can be assimilated to  time. 
On an energy surface $\E_{\hat{H} = \hat{H}_0}$, consider the Poincar{\'e}
section $s=0$  parametrized by $(\theta,r) \in {\T} \times \R$
$$u := \hat{H}_0 - { r^2 \over 2} - \varepsilon H(\theta,r,0).$$

It was proved by Robinson \cite{robinson} that for the generic $H \in C^{\infty} ( {\T} \times \R \times
{\T}, \R) $, the flow $X_{ \hat{H}}^t$ satisfies the Kupka Smale
condition. In this case, for $\varepsilon$ small enough, there exists
$R > 0$ such that the projection on $(\theta,r)$, $|r| \leq R$, of the Poincar{\'e}
section map is a monotone twist map, denote it  $f$. From KAM
theory, we know that there exists an
invariant annulus for this map inside $\T \times [-{R \over 2}, {R
  \over 2}]$; denote this annulus by $A_0$. 

\subsection{} \label{mather} A compact set $\A$  in the annulus $A_0 \subset \T \times \R$ is said to be $f$-ordered if it is invariant by $f$, if it projects injectively on ${\T}$ and if the restriction of $f$ to $\A$ preserves the natural order given by the projection.  From the Aubry-Mather theory it follows from the fact that $f$ is a monotone twist map that for every $\alpha \in [\a_1, \a_2] $, where $\a_i$, $i=1,2$ are the rotation numbers on the boundaries of $A_0$, there exists a minimal $f$-ordered set $\A_{\alpha}
\subset A_0$. If $\a$ is irrational  the set $\A_{\alpha}$ is either a Cantor set or a continuous graph $(\theta, \psi(\theta))$. In both cases order preserving implies that the restriction $f_\a$ of $f$ to $\A_\a$ is   semi-conjugate to the irrational rotation $R_\a$ just like orientation preserving homeomorphisms of the circle with irrational rotation number are.  The case of a Cantor set corresponding to a homeomorphism called a {\sl Denjoy counter-example} that displays a wandering interval, i.e., an interval disjoint from all its iterates. 

We will be chiefly interested in this paper with the case where $\A_\a$ is a Cantor set. The semi-conjugacy between $f_\a$ and the rotation is obtained in the following way: we choose a point $u$ on the Cantor set and we project it to a point on the circle, say $0$. Then to each iterate of $u$ by $f_\a$ we associate the iterate of $0$ by $R_\a$. Then we extend by continuity using minimality and bridging between any two points of the Cantor set separated by a gap by projecting them on the same point. Doing so we get a map $h$ that is continuous from $\A_\a$ onto $\T$ and satisfies $h \circ f_\a = R_\a \circ h$. Besides the map $h$ is injective except on endpoints of gaps, hence we can define an inverse $h^{-1}$ except at the points of countably many orbits of the rotation corresponding to different orbits of gaps on $\A_\a$; we will call each such orbit a {\sl hole}. At the points where $h$ is not invertible, we take $h^{-1}$ to be the right extremity of the corresponding gap. The map $h^{-1}$ thus defined will be right semi-continuous. 

We recall that $\A_\a$ lies on a Lipschitz graph and that the projection on $\T$ is bi-Lipschitz. This will be of importance in \S \ref{redd} as we will  reduce the study of special flows over Aubry-Mather sets to the study of special flows over rotations on the circle.


Finally, to any Aubry-Mather set $\A_{\a}$ of the
Poincar{\'e} section map defined above corresponds for the flow $X_{
  \hat{H}}^t$ a minimal invariant set ${\s}_{\a}$ on which the dynamics
is  a suspension over
  $f_{\a}$ with a constant ceiling function equal to one (since $ds / dt =1$). The only
  invariant probability measure on  ${\s}_{\a}$ is the
  normalized product
  of the invariant measure of $f_\a$ on $\A_\a$ with the Lebesgue measure on the
  fibers. Hence the spectrum of  the flow $X_{
  \hat{H}}^t$ on $\s_\a$ is not continuous since the function $e^{i 2 \pi s}$ is an eigenfunction 
of the flow. The eigenfunction being continuous the flow is not topologically weak mixing (Cf. \S \ref{weak}). 

\vspace{0.15cm}

\subsection {\sl Reparametrizations.} \label{rep} We address the following question: {\sl How sensitive is the dynamics on ${\s}_{\a}$ under time change of the flow, i.e., if
all the orbits of $X_{ \hat{H}}^t$ are kept intact but the
 velocity along them is modified?} 

A smooth reparametrization of a flow given by a smooth vector field
is done by multiplication of  the original vector field by some smooth and strictly positive real function $\phi$. In the case of a Hamiltonian system, if we restrict our study to an energy surface $\E_{\hat{H} = \hat{H}_0}$ of the flow  $X_{ \hat{H}}^t$,
 a time change can be performed using the
Hamiltonian 
\begin{eqnarray} \label{amh}
\hat{H}_{\phi} := \phi(\theta,r,s,u) (\hat{H}(\theta,r,s,u) -\hat{H}_0) ,
\end{eqnarray}
for some smooth real function $\phi$, $\phi > 0$.  The energy
surface $\E_{\hat{H}_0}$ is indeed invariant by the Hamiltonian flow
corresponding to $\hat{H}_{\phi}$ denoted by $X_{ {\hat{H}}_{\phi}}^t$; and if
we write the vector field of $X_{ {\hat{H}}_{\phi}}^t$  on 
$\E_{\hat{H}_0}$, we notice that it is exactly the vector field of
$X_{ {\hat{H}}}^t$ multiplied by $\phi$.

It is a general fact that the time change of a uniquely ergodic system is uniquely ergodic (Cf. for example \cite{parry}).
Here the unique invariant probability measure of the reparametrized flow is absolutely continuous with respect to the original measure and has density ${1 \over \phi}$.

In the sequel, we will discuss whether a reparametrization given by (\ref{amh}) can yield a continuous spectrum 
 on the laminations ${\s}_{\alpha} \subset \E_{\hat{H}_0} $. Since we are considering perturbations of completely integrable Hamiltonians we will pay a special attention to  reparametrizations that are small perturbations, that is, in time change functions $\phi$ that are close to one.

\subsection{ \sl Continuous spectrum.}  \label{weak} 
We say that a measure preserving flow 

\noindent $( T^t, X, \mu )$ has a {\sl continuous spectrum}  if and only if it does
not have an
eigenfunction,
 i.e., a measurable non constant complex function $h$ such that 
 $ h \left( T^t x \right) =   e^{i \lambda t} h \left( x
\right),$ at almost every $x \in X$ for some $\lambda \in \R$.

 An equivalent property is {\sl weak mixing}: We say that a measure preserving flow $( T^t, X, \mu )$ has the weak mixing property  if and only if 
 for all measurable sets $A$ and $B$ 
$$ \mu ( T^{-t} A \bigcap B ) \longrightarrow \mu(A) \mu (B), $$
when $|t|$ goes to infinity on a set of  density one over $\R$.

 For the equivalence between the definitions see for example the book by Parry \cite{parry}.

A flow is said to be {\sl topologically weak mixing} if it does not have a continuous non constant eigenfunction.
\vspace{0.15cm}

\subsection{\sl The case of invariant tori.} \label{tores} In the case  where the   invariant set $\A_\a$ of the Poincar\'e map is a smooth circle and $f_\alpha$ is $C^{\infty}$
conjugate to the translation $R_\alpha$ (KAM invariant tori), the question of possibly having a continuous spectrum after time change was asked by Kolmogorov  and treated by himself \cite{Kkk}, Shklover \cite{shklover}, Herman \cite{He1} and others. In \cite{weakmixing} the following dichotomy was proved:

\vspace{0.2cm}

\noindent {\sc Theorem.}  {\sl Let $R_{\a}$ be a minimal translation on  ${\T}^d$. Then
  either one
  of two possibilities hold: 

\noindent $(i)$ The vector $\a$ is  Diophantine, and for any $\phi \in C^{\infty}({\T}^{d+1}, \R_{+}^{*})$, the reparametrization  of
the flow $R_{t(1,{\a})}$, with speed $ \phi$  is
$C^{\infty}$  conjugate to a translation flow on $ {\T}^{d+1}.$ 

\noindent $(ii)$ The vector $\a$ is Liouville (i.e. not Diophantine), and for a dense $G_{\delta}$
of $\phi \in  C^{\infty}({\T}^{d+1}, \R_{+}^{*})$, 
 the reparametrization  of
the flow $R_{t(1,{\a})}$, with speed ${ \phi}$  is weak mixing
  (for its unique invariant measure).}

\vspace{0.2cm}

Alternative (i) was proven by Kolmogorov for $d=1$ and generalized to any dimension by Herman. 

\subsection{ \sl The case of laminations.} \label{trous}


\begin{defi} \label{PPPP} 
{ If the Aubry-Mather set has exactly  $k$ orbits of wandering intervals, we
say it is an} {Aubry-Mather set with $k$ holes}. 
\end{defi}

 {
Assume  $I$ and $J$ are two gaps whose endpoints lie on distinct orbits of ${f_\a}$ on ${| \A_\a}$
and let $\beta_I$ and $\beta_J$ be their point projections on the circle by some semi-conjugacy  of $f_\a$ to $R_\a$.
It is easy to see that the difference $\beta_{I} - \beta_{J}$ does
not depend on the choice of the semi-conjugacy.}

\begin{defi} { Holes in a general position.} \label{amhhhh}
{ We
say that the { holes} of $\A_\a$, corresponding to $I$ and $J$}, { are}  in a general position  { if } $\beta = \beta_I - \beta_J$ { satisfies}
\begin{eqnarray} \label{sstar} \hspace{4cm} ||q_n \beta|| \mathop{ \not\longrightarrow} \limits_{n \rightarrow \infty} 0, \hspace{6cm}  \end{eqnarray}
{where $q_n$, $n \in \N$ are the denominators in the successive convergents of $\a$, and $||.||$ is the distance to the closest integer.   
} \end{defi}

We remind that it is possible to define the sequence $q_n$ in the following unique way:  $q_{-1}= q_0 =1$, and for every $n \in \N$  
\begin{eqnarray} \label{amapproximations} \hspace{2cm} \parallel q_n \a \parallel
  < \parallel k \a \parallel, \hspace{1.5cm} \forall k < q_{n+1}. \end{eqnarray}
The sequence $q_n$ can equally be defined as satisfying the recurrence relation
\begin{eqnarray} q_n = a_n q_{n-1} + q_{n-2} \ \ {\rm for} \ \ n \geq 1, \ \ q_{-1} = q_0=1,  \label{amqn}
\end{eqnarray} 
where the $a_n$ are the  partial quotients in the continued fraction expansion of $\a$ (Cf. for example \cite{He2}, Chapter V).

\vspace{0.2cm} 

\noindent {\sc Remark.} 
{ Given an irrational number $\a$,  the set of numbers} $\beta$,  
{  satisfying} $(\ref{sstar})$
{ is a dense $G_{\delta}$ in $\R$ of full measure on any compact interval.}  


\vspace{0.2cm}

\subsection{\sl Hyperbolic Aubry-Mather sets.}  Our first result,
related to Aubry-Mather sets with only one hole, will be stated just
for {\sl hyperbolic} Aubry-Mather sets (i.e. admitting stable and unstable bundles) and will be based on the following facts: 

\vspace{0.2cm} 

\noindent 2.7.1. The assumption of hyperbolicity  is relevant  since
P. Le Calvez has shown that under the same
generic condition that  $X_{ \hat{H}}^t$ should be Kupka Smale,
it is true for an open and dense set of $\a \in \R$
that the Aubry-Mather set with rotation number $\a$ is hyperbolic
\cite{Patrice}. 

\vspace{0.2cm} 

\noindent 2.7.2.  When
this is the case, A. Fathi
shows that the Hausdorff dimension of the Aubry-Mather set is
equal to zero \cite{Fathi} (he shows even more, that the union of the
hyperbolic Aubry-Mather sets has Hausdorff dimension zero). 

\vspace{0.2cm}

\noindent 2.7.3. The sizes of  $f_{\a}^i (I)$ decrease to zero since they are disjoint intervals on a Lipschitz graph. The hyperbolicity assumption
on $\A_{\a}$ forces these sizes to decrease in fact geometrically.


\section{The results.} \label{amresults}

\subsection{} The notations are as in the  precedent section: the flow corresponding to the Hamiltonian given in (\ref{hamo}) is denoted by   $X_{{\hat{H}}}^t$. 
 The flow arising from the reparametrization (\ref{amh}) of $X_{{\hat{H}}}^t$, on the energy level $\E_{\hat{H} = \hat{H}_0}$ and with speed $\phi$, is denoted by   $X_{{\hat{H}_{\phi}}}^t$. On the energy surface $\E_{\hat{H} = \hat{H}_0}$ we consider the  Poincar\'e section $s=0$ and an annulus $A_0$ on this section invariant by the Poincar\'e return map (\S \ref{set}). An  Aubry-Mather set corresponding to a frequency $\a$ is denoted by $\A_\a$, and $\s_\a$ denotes the  lamination  over  $\A_\a$ invariant by the flow $X_{{\hat{H}}}^t$  as well as by the flow  $X_{{\hat{H}_{\phi}}}^t$ . 
 Finally, the time change function $\phi$, defined on $\T \times \R \times \T \times \R$ is assumed to be of class at least $C^1$ and strictly positive. Naturally, the reparametrized flow is as close to the initial flow as $\phi$ is close to 1.

\subsection{} \label{amoneone} Recall Definition \ref{PPPP} of (\ref{trous}) on the number of holes of an Aubry-Mather set. 


\vspace{0.2cm}

\noindent {\sc Theorem.}  {\sl If $\A_{\alpha}$ is a hyperbolic Aubry-Mather set and has  one hole, then, for any time change
function $\phi$ of class $C^1$ with mean value 1, the restriction of the flow 
  $X_{{\hat{H}}_{\phi}}^t$
to ${\s}_{\alpha} $ is $C^0$ conjugate to the initial flow, i.e. to a
suspension flow above $(\A_\a, f_\alpha)$ with a constant suspension function
equal to one.} 

\vspace{0.2cm}

It immediately follows that, under these hypothesis on $\A_\a$, the restriction  of the reparametrized flow to  $\s_\a$ is never topologically weak mixing.


\subsection{} \label{ammanymany} In the case $\A_\a$ is an Aubry-Mather sets with more than one hole ($\A_\a$ hyperbolic or not),  
we have with the notion of general position given in Definition \ref{amhhhh} of (\ref{trous})

\vspace{0.2cm}

\noindent {\sc Theorem.}  {\sl 
If $\A_{\alpha}$ is an Aubry-Mather set with at
least two holes in a general position, then arbitrarily close (in the
$C^{\infty}$ topology) to the
constant function equal to one, there exists a function
$\phi \in C^\infty ({(\T \times \R )}^2, \R_+^*)$ such that the  restriction of the reparametrized flow  $X_{{\hat{H}}_{\phi}}^t$
to ${\s}_{\alpha}$ is weak mixing for the unique invariant measure.}
 
\vspace{0.2cm}

\subsection{} \label{co} Assume now that $\A_{\a}$ is any given Aubry-Mather set on the Poincar\'e section. We can always assume that the semi-conjugacy to $R_\a$  in \S  \ref{mather} projects an orbit of a hole of $\A_\alpha$ by $f_\a$ to the orbit of 0 by $R_\a$.  If we consider a two fold covering of the Poincar\'e  section, we obtain an Aubry-Mather set $\overline{\A}_\a$ for the Poincar\'e map of the Hamiltonian flow given by $\overline{H}
(\theta,r,s,u) := \hat{H} (2\theta,r,s,u)$. The set  $\overline{\A}_\a$ now has  at least two holes with one projecting  in 0 and one in ${1 \over 2}$. From the last theorem we will be able to deduce the following 

\vspace{0.2cm}

\vspace{0.2cm} 

\noindent {\sc Corollary.} {\sl
  Arbitrarily close (in the
$C^{\infty}$ topology) to the
constant function equal to one, there exists a $C^{\infty}$ function
$\phi$ such that the  restriction of the flow  $X_{ {\overline{H}}_{\phi}}^t$
to the lamination $\overline{\s}_{\alpha}$ over  $\overline{\A}_\a$ is weak mixing for the unique invariant measure.}

\vspace{0.2cm}

\section{Proofs.} \label{amproofs}

\subsection{\sl Reduction to special flows.}  On the energy surface     
 $\E_{\hat{H} = \hat{H}_0}$, the flow $X^t_{\hat{H}}$ is a constant time suspension over the Poincar\'e map $f$ of the Poincar\'e section $s=0$. We restrict $f$ to the invariant annulus $A_0$ on which we introduce a pair of coordinates $(\tilde{r}, \tilde{\theta})$.  The reparametrized flow 
 $X_{{\hat{H}}_{\phi}}^t$ can be viewed as a special flow over $(f, A_0)$ with a ceiling function $\psi$ given by 
\begin{eqnarray} \label{speciale} \psi(\tilde{r}, \tilde{\theta}) := \int_0^1 {1 \over \phi(X^s_{\hat{H}}(\tilde{r}, \tilde{\theta},s))} ds. \end{eqnarray}
 In the proof of the Theorems, it will be more convenient to work with special flows rather than with reparametrizations.  For Theorem \ref{ammanymany}, we will construct a special flow above $(f_\a,A_\a)$ that is weak mixing and show that the special function  we will use can be obtained from a smooth time change function $\phi$ via the formula (\ref{speciale}).

\subsection{\sl Proof of Theorem \ref{amoneone}.} \label{amone} By the formula (\ref{speciale}), if $\phi$ is of class $C^1$ then so will be the function $\psi$. Hence, the proof of the theorem will be accomplished if we
 show that the special  flow over $f_\a$ and under any function that is the restriction over $\A_\a$ of a function $\psi \in C^1(A_0, \R_+)$  is $C^0$ conjugate to a suspension
 with a constant ceiling function. This in its turn will follow if we prove the existence of a continuous solution on $\A_\a$ to the equation 
\vspace{0.2cm}

 \begin{eqnarray} \label{(E)} \xi - \xi \circ f_\a  = \psi - \int_{\A_\a} \psi(\tilde{r}, \tilde{\theta}) d \mu_\a, \end{eqnarray}
where $\mu_\a$ is the unique invariant measure on $\A_\a$ by $f_\a$ (Cf. for example \cite{HK} Chapter 4).   
\vspace{0.2cm}

 By the Theorem of Gottschalk and Hedlund (\cite{gott} or
\cite{He2} Chapter IV), equation (\ref{(E)}) will have a continuous solution if we
prove that $|\psi_m(x) - \psi_m(x')|$ is uniformly bounded for $m \in
\N$, $x,x' \in \A_{\a}$ (here $\psi_m(x)$ denotes the Birkhoff sums of
$\psi$ relative to $f_\a$, i.e. $\psi_m(x) = \psi(x) + \psi ( f_\a x) + ... + \psi(f_\a^{m-1} x)$.

\subsection{\sl Reducing to special flows over rotations of the circle.}  \label{redd} Assume now that $\A_{\a}$ is  a hyperbolic Aubry-Mather set with only one hole. To the special flow over $(\A_\a,f_\a)$ and under the function $\psi$ we associate a special flow over $({\T}^1, R_\a)$ and under a function  $\varphi$ 
and showing that the sums $S_m \psi(x) - S_m \psi(x') $ above $f_\a$ are uniformly bounded is equivalent to showing that $S_m \varphi(\theta) - S_m \varphi(\theta') $ above $R_\a$ are uniformly bounded. The function $\varphi$  has the following properties

\vspace{0.2cm} 

\noindent {\sl i)}  The function $\varphi$ is strictly positive and for any  $0 \leq a < b < 1 $
\begin{eqnarray} \label{amamam11} \varphi(b) -\varphi(a)=
\sum_{R^k_\a (0) \in ]a,b]} \Delta_k, \end{eqnarray}
where the number $\Delta_k$ is the difference between the values of $\psi$ at the right and left endpoints of the $k^{\rm th}$ gap. The reason for (\ref{amamam11}) is that we assumed  that $(f_\a, \A_\a)$ is hyperbolic: which implies (Cf. \S 2.6.2) that the union of the gaps on the Lipschitz graph where $\A_\a$ lies has full Lebesgue measure. Since $\psi$ is of class $C^1$ we have that the variations of $\psi$ are concentrated on the gaps and since the projection on $\T$  is bi-Lipschitz we get that the corresponding variations of $\varphi$ are concentrated on the orbit of $0$ as stated in (\ref{amamam11}).

Since we assumed that $(f_\a, \A_\a)$ is hyperbolic we also have (Cf. \S 2.6.3.) that the size of the $k^{\rm th}$ gap decreases geometrically as $k \rightarrow \pm \infty$.  The function $\psi$ being Lipschitz we obtain the following

\vspace{0.2cm} 

\noindent {\sl ii)}   The sequence $ {\lbrace |{\Delta}_k| \rbrace}_{k \in \Z}$  decrease
geometrically when $|k| \rightarrow \infty$, e.g. there exist
 $ C > 0$ and $ 0< \Delta < 1 $ such that for any $k \in \Z$
\begin{eqnarray}  
|{\Delta}_k| \leq C {\Delta}^{|k|}. \end{eqnarray}
We will actually need the following weaker property on the $\Delta_k$:
\begin{eqnarray}  
\sum_{k \in \Z} |k {\Delta}_k| < + \infty .  \label{amamam33} \end{eqnarray}

\vspace{0.2cm} 

\noindent {\sl iii)}  Moreover, we clearly have
\begin{eqnarray} \label{amamam22} \sum_{i= -\infty}^{+\infty} \Delta_i = 0. \end{eqnarray}

 In conclusion, to obtain the theorem we just need to prove the following
\vspace{0.2cm} 

\noindent {\sc Proposition.}  {\sl Given any rotation on the circle $R_\a$  of the circle and a real function $\varphi$ satisfying 
  {\sl i) --- iii)}, we have 
$$\sup_{m\in {\N}^*} \sup_{(\theta,\theta')\in \T} |\sum_{i=0}^{m-1} \varphi(\theta'+i\a) -
\sum_{i=0}^{m-1} \varphi(\theta+ i \a) | < + \infty.$$
}

\vspace{0.2cm} 

\noindent{\sl Proof. } For the proof of the proposition we can assume
 that the integral of $\varphi$ vanishes. From (\ref{amamam11}) we have that the derivative of $\varphi$ in
the sense of the distributions is
$$ D \varphi = \sum_{k \in \Z} \Delta_k \d_{k\a},$$
where $\delta_z$ denotes the Dirac measure concentrated on $z$.

Define for every $k \in \Z$
$$\sigma_k~:= \sum_{j= -\infty}^{k} \Delta_j = - \sum_{j= k+1}^{+\infty} \Delta_j$$
and denote by $e_k$ the function on the circle, of zero integral satisfying 
$$D e_k = \d_{k\a} - \delta_{(k+1) \a}.$$

Since $\Delta_k = \sigma_{k} - \sigma_{k-1}$, we have
$$\varphi = \sum_{k \in \Z} \sigma_k e_k.$$

Now, if we use the usual notation $S_me_{k}$ for the Birkhoff sums of
$e_k$ above $R_\a$ we notice that for any $m >0$ the following is true
$$D S_me_{k} = \d_{k-(m-1)\a} - \d_{(k+1)\a}$$
hence
$$ {|| S_me_{k} ||}_{L^{\infty}} \leq 1$$
which implies
\begin{eqnarray*} {|| S_m \varphi ||}_{L^{\infty}} &\leq& \sum_{k
    \in \Z} | \sigma_k | \\  
&\leq& \sum_{i \in \Z} (1 + |i|) |\Delta_i|. \end{eqnarray*}

Condition (\ref{amamam33}) concludes the proof of the proposition and hence of Theorem \ref{amoneone}. \carre

\subsection{\sl Proof of Theorem \ref{ammanymany}.}  The proof of weak
mixing we will give in this section is similar to the one produced by 
Katok and Stepin in \cite{KS} for interval exchange transformations.

\subsection{\sl Reduction to special flows.}  \label{amam2} Assume $\A_{\a}$ is an Aubry-Mather set having two orbits of wandering
intervals one in 0 and the other in $\beta$, and that $\beta$
satisfies condition $(\ref{sstar})$ of Definition \ref{amhhhh}. Define the set 
 $$ \X_{\a} = \lbrace \epsilon  \in ]0,1[, \ \ {1 \over   \epsilon} \not\in \Q + \a \Q \rbrace.$$ 
Theorem \ref{ammanymany} will follow if we prove the following

\vspace{0.2cm} 

\noindent {\sc Proposition} {\sl
 For any $\epsilon \in  \X_{\a} $, 
the special  flow over $R_{\a}$
  with the ceiling function  
\begin{eqnarray} {\chi}^{\epsilon} = (1- \epsilon) {\chi}_{[0, \beta[}
+ (1+ \epsilon ) {\chi}_{ [ \beta, 1 [ } \label{qui} \end{eqnarray}
is weak mixing (for the unique
  invariant measure).}
 
\vspace{0.2cm}

Indeed, the Aubry-Mather set $\A_\a$ has a gap in $0$ and a gap in $\beta$ that separate it in two parts. Hence we can find two disjoint open sets of ${(\T \times \R)}^2 $ each one containing a part of the lamination $S_\a$. We then choose the time change function $\phi^{\epsilon}$ in (\ref{amh}) of class $C^\infty$ and constant on each one of the two open sets with values ${(1-\epsilon)}^{-1}$ on one, and ${(1+\epsilon)}^{-1}$ on the other. By the formula (\ref{speciale}), the restriction of the reparametrized flow to $S_\a$ can be viewed as the special flow over $f_\a$ with the ceiling function (\ref{qui}).  By semi-conjugacy, the weak mixing of this special flow will follow from Proposition \ref{amam2}. Furthermore, it is clearly possible to construct the $\phi^\epsilon$ as close to 1 in the $C^\infty$ topology as $\epsilon$ goes to zero in (\ref{qui}).

\vspace{0.2cm}

\noindent {\sl Proof of Proposition \ref{amam2}.} We will use a classical general lemma on weak mixing for 
special flows, the proof of which can be found for example in \cite{CFS}.
 In the lemma  $\lbrace T^t \rbrace $ can be any 
special flow constructed from an ergodic automorphism $(T,L,\mu)$ of a
Lebesgue space $M$ and under  a summable function $f > 0$. For the special flow we will consider the normalized measure $ {1 \over \int f d\mu} d\mu ds$ where $ds$ denotes Lebesgue measure on the fibers.

\vspace{0.2cm} 

\noindent {\sc Lemma.} {\sl The flow  $\lbrace T^t \rbrace$ has a continuous spectrum if and
only if for any $\lambda$ in ${\R}^{*}$ the equation
\begin{eqnarray} \label{eee} h(T(x)) = e^{i \lambda f(x)} h(x),
\end{eqnarray}
does not have a non zero measurable solution $h$. }

\vspace{0.2cm}

 In our case (\ref{eee}) becomes 
\begin{eqnarray} \label{eeee} h(x+\a) = e^{i \lambda \chi^\epsilon (x)} h(x), \end{eqnarray}
where 
 \begin{eqnarray*}
e^{i \lambda \chi^\epsilon (x)}&=&  e^{i \lambda (1- \epsilon) } \ \ \
if x \in [0, \beta[, \\
e^{i \lambda \chi^\epsilon (x)}&=&  e^{i \lambda (1+ \epsilon) } \ \ \
if x \in [ \beta,1[.
\end{eqnarray*}

 We will use the following

\vspace{0.2cm} 

\noindent {\sc Lemma.} {\sl Let $\a$ be an irrational number.  If the number  $\beta \in ]0,1[$ satisfies $(\ref{sstar})$, and $w(x)$ is a  complex
function defined on ${\T}^1$ such that 
\begin{eqnarray*}
w(x) &=& z_1, \ \ \ if x \in [0, \beta[, \\
w(x) &=& z_2, \ \ \ if x \in [\beta,1[, 
\end{eqnarray*}
with $|z_1|=|z_2|=1$ and $z_1 \neq z_2$; then the equation 
\begin{eqnarray} \label{amhomologie} 
h(x + \a) = w(x) h(x), \end{eqnarray}
does not admit any measurable solution $h$.}

\vspace{0.2cm}

This Lemma was stated and proved by Katok and Stepin in \cite{KS} in the case where $\a$ is not of constant type (the sequence $a_n$ in \ref{amqn} is unbounded). Their proof was based on fast cyclic approximations for irrational rotations of non-constant type. The proof of the absence of solutions in the case where $\a$ is of constant type is due to Veech \cite{Veech} (a complete discussion with proofs can be found in \cite{merrill}).

The last lemma implies that a necessary condition on an eventual eigenvalue
$\lambda$ for the special flow of Proposition \ref{amam2}  is that
$$ \lambda \epsilon = l \pi, \  {\rm for \  some \  } l \in \Z.$$
But then a corresponding solution of (\ref{eeee}) would satisfy
$$h(x + \a) =   e^{i {l \pi }}   e^{i {l \pi \over \epsilon} } h(x),$$
hence  $ e^{i {l \pi }}   e^{i {l \pi \over \epsilon} }$ should be an eigenvalue of $R_\a$, that is 
$$ 
 l + {l \over \epsilon} = k \a + 2p,$$
for some integers $k$ and $p$, which contradicts the fact that $\epsilon \in \X_\a$.

Therefore equation (\ref{eeee}) does not have non trivial solutions and the
special flow of Proposition \ref{amam2} is weak mixing. Theorem \ref{ammanymany} is thus proved. \carre

\subsection{\sl Proof of the Corollary \ref{co}.} \label{am3} 
The corollary of Theorem \ref{ammanymany} will follow if we prove the following

\vspace{0.2cm} 

\noindent {\sc Proposition.} {\sl
 For any $\a \in \R - \Q$, and any $\epsilon \in \X_\a$,  the special  flow over $R_{\a}$
  with the ceiling function  $\chi^\epsilon = (1- \epsilon) \chi_{[0,
  {1 \over 2}[}
+ (1+ \epsilon ) \chi_{[ {1 \over 2},1[} $ is weak mixing for its unique
  invariant measure.} 

\vspace{0.2cm}

\noindent {\sl Proof. } From Proposition \ref{amam2}, we just have to prove that 0 and ${1 \over 2}$ are independent with respect to $\a \in \R-\Q$.

Since two consecutive denominators of the convergents of $\a$ are
relatively prime (this can be easily obtained by recurrence from
(\ref{amqn})), one at least among two consecutive $ q_{n}$ and $ q_{n+1}$ is odd.
We  extract by this means a subsequence of the best approximations of
$\a$ for which ($\ref{sstar}$) of Definition \ref{amhhhh} holds for ${1 \over 2}$. Proposition \ref{am3} is thus
proved. \carre

\vspace{0.8cm}

\noindent{\sc Acknowledgment.} I wish to thank 
Patrice Le Calvez, Fran{\c c}ois Parreau and J-C Yoccoz for 
valuable remarks and simplifications.

\vspace{0.1cm} 

I am grateful for the late Michael Herman who suggested the problem to me.

\vspace{0.8cm}

\bibliographystyle{plain}
\bibliography{these}

\vspace{0.4cm} 

\noindent Bassam Fayad

\noindent LAGA,  Universit\'e Paris 13

\noindent E-mail: fayadb@math.univ-paris13.fr

\end{document}